\documentclass[iicol]{sn-jnl}

\jyear{2021}%

\theoremstyle{thmstyleone}%
\newtheorem{theorem}{Theorem}
\theoremstyle{thmstyletwo}%
\newtheorem{remark}{Remark}%

\theoremstyle{thmstylethree}%

\def\kappagamma{\small \bigl[\!{\begin{array}{c}\kappa\\[-1ex] \gamma\end{array}}\!\bigr]} 
\def\xizeta{\footnotesize\bigl[\!{\begin{array}{c}\xi\\[-0.2ex] \zeta\end{array}}\!\bigr]} 
\def\kappatil{\underline{\widetilde{\kappa}}}
\def\gammatil{\underline{\widetilde{\gamma}}}
\def\indset{I}
\def\freqset{U}
\def\phiTay{\check{\phi}}
\def\psiTay{\check{\psi}}
\def\multip{m}
\usepackage{amsmath,amssymb,comment}
\newcommand{\inarxiv}[1]{#1}

\begin{document}

\title[vibro-acoustography]{On the inverse problem of vibro-acoustography}


\author*[1]{\fnm{Barbara} \sur{Kaltenbacher}}\email{barbara.kaltenbacher@aau.at}







\abstract{The aim of this paper is to put the problem of vibroacoustic imaging into the mathematical framework of inverse problems (more precisely, coefficient identification in PDEs) and regularization. We present a model in frequency domain, prove uniqueness of recovery of the spatially varying nonlinearity parameter from measurements of the acoustic pressure at multiple frequencies, and derive Newton as well as gradient based reconstruction methods.   
}

\keywords{vibro-acoustic imaging, inverse problem, coefficient identification, regularization}


\pacs[MSC Classification]{35R30,65J20}

\maketitle


\section{Introduction}
Vibro-acoustography by means of ultrasound  was developed \cite{FatemiGreenleaf1998,FatemiGreenleaf1999} to achieve
the high resolution by high frequency waves while avoiding the drawbacks of scattering from small inclusions and of stronger attenuation at higher frequencies.
\begin{figure}
\includegraphics[width=0.48\textwidth]{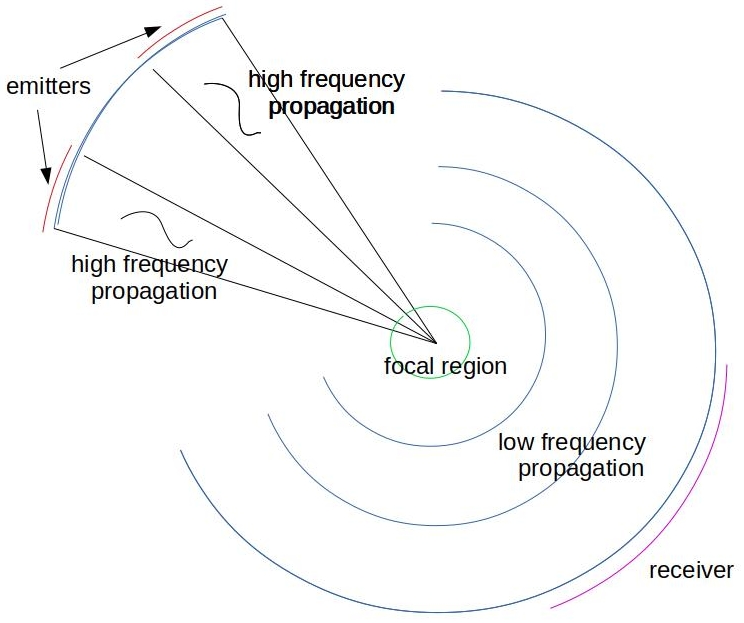}
\caption{schematic of the experimental setup\label{fig:experiment}}
\end{figure}
The experiment for image acquisition is illustrated in Figure \ref{fig:experiment}: 
Two ultrasound beams of high and slightly different frequencies $\omega_1$ and $\omega_2$ are excited at two parts $\Sigma_1$, $\Sigma_2$ of an array of piezoelectric  transducers (emitters). They interact nonlinearly at a focus and this interaction excites a wave that basically propagates at the difference frequency $\omega_1-\omega_2$ and is eventually measured by a receiver array $\Gamma$ (in experiments consisting of hydrophones, in imaging this would be a piezoelectric transducer array as well).
After each measurement, the focal region is shifted to scan the overall region of interest.
Inhomogeneity of the medium leads to spatial dependence of two coeffients in the governing models: The speed of sound $c=c(x)$ and the nonlinearity parameter $\gamma=\gamma(x)$.
Both parameters are susceptible to local variations in the acoustic medium (e.g., human tissue in medical applications) and thus their reconstruction yields a spatial image of the region of interest.
In case of reconstructing the individual coefficients $c=c(x)$ or $\gamma=\gamma(x)$, this is related to ultrasound tomography and nonlinearity parameter imaging, respectively, cf. e.g., \cite{Bjorno1986,Cain1986,IchidaSatoLinzer1983,SPIE2021} and the citing literature.

A modeling and simulation framework for this methodology has been devised in 
\cite{Malcolmetal2007,Malcolmetal2008}. 
In this paper we put an emphasis on the inverse problem of reconstructing $c=c(x)$ and $\gamma=\gamma(x)$. 

\section{Model}
Two ultrasound beams with acoustic velocity potentials $\phi_1$, $\phi_2$ are excited by transducers 
and their interaction in turn excites a wave field with velocity potential $\psi$. This is discribed by a system of PDEs with inhomogeneous Neumann conditions
\begin{eqnarray}
&&\partial_t^2 \phi_k-c^2\Delta \phi_k=0 \mbox{ in }\Omega,\nonumber\\ 
&&\partial_\nu \phi_k=g_k \mbox{ on } \Sigma_k\,,\ k\in\{1,2\}\label{eq:phik}\\
&&\partial_t^2 \psi-c^2\Delta \psi= \tilde{f}(\phi_1,\phi_2,\gamma,c) \mbox{ in }\Omega
\label{eq:psi}\\
&&\tilde{f}(\phi_1,\phi_2,\gamma,c)=\nonumber\\
&&\partial_t\left(\left\vert\nabla(\phi_1+\phi_2)\right\vert^2 + \frac{\gamma-1}{2c^2}
\left\vert\partial_t(\phi_1+\phi_2)\right\vert^2\right)
\nonumber
\end{eqnarray}
see \cite{Westervelt:63} for the derivation of the nonlinear forcing $f$.
In here, $c=c(x)$ and $\gamma=\gamma(x)$ are the spatially varying sound speed and nonlinearity parameter, respectively, and the manifold $\Sigma_k$ represents the emitting transducer array with given time harmonic excitation $g_k(x,t)=\hat{g}_k(x)\, e^{\imath\omega_k t}$.

The system \eqref{eq:phik}, \eqref{eq:psi} is not fully nonlinear but the task of its solution can be decoupled into two linear subproblems: First compute $\phi_1,\phi_2$ from \eqref{eq:phik}, then insert them into the right hand side of \eqref{eq:psi}, and finally solve \eqref{eq:psi} for $\psi$.

\paragraph{Transformation into frequency domain}
Linearity of the subproblems allows to easily transfer the time domain formulation \eqref{eq:phik}, \eqref{eq:psi} into frequency domain.
With the time harmonic ansatz $\phi_k(x,t)=\hat{\phi}_k(x)\, e^{\imath\omega_k t}$, $\psi(x,t)=\Re\Bigl(\hat{\psi}(x)\, e^{\imath(\omega_1-\omega_2)t}\Bigr)$, where the latter is induced by real-valuedness of the right hand side of \eqref{eq:psi}
\begin{eqnarray*}
&&f(\phi_1,\phi_2,\gamma,c)\\
&&=\partial_t\Bigl[\left\vert\nabla\hat{\phi}_1\right\vert^2+\left\vert\nabla\hat{\phi}_2\right\vert^2
+2\Re\left(\nabla\hat{\phi}_1\cdot\nabla\overline{\hat{\phi}_2}\,  e^{\imath(\omega_1-\omega_2)t}\right)\\
&&\qquad+\frac{\gamma(x)-1}{2c(x)^2}\Bigl(\left\vert\omega_1\hat{\phi}_1\right\vert^2+\left\vert\omega_2\hat{\phi}_2\right\vert^2\\
&&\qquad\qquad+2\omega_1\omega_2\Re\left(\hat{\phi}_1\overline{\hat{\phi}_2}\, e^{\imath(\omega_1-\omega_2)t}\right)
\Bigr)\Bigr]\\
&&=2(\omega_1-\omega_2)\\
&&\ \cdot\Re\Bigl(\imath\left(
\nabla\hat{\phi}_1\cdot\overline{\nabla\hat{\phi}_2} 
+\omega_1\omega_2 \frac{\gamma(x)-1}{2c(x)^2}\hat{\phi}_1\overline{\hat{\phi}_2}\right) e^{\imath(\omega_1-\omega_2)t}\Bigr),
\end{eqnarray*}
we get 
\begin{eqnarray}
&&-\frac{\omega_k^2}{c(x)^2} \hat{\phi}_k-\Delta \hat{\phi}_k=0 \mbox{ in }\Omega \nonumber\\
&&\partial_\nu \hat{\phi}_k=\hat{g}_k \mbox{ on } \Sigma_k\,,\quad k\in\{1,2\}
\label{eq:phikhat0}\\
&&
-\frac{(\omega_1-\omega_2)^2}{c(x)^2} \hat{\psi}-\Delta \hat{\psi}=
f(\hat{\phi}_1,\hat{\phi}_2,\gamma,c)\mbox{ in }\Omega\quad\label{eq:psihat0}\\
&&f(\hat{\phi}_1,\hat{\phi}_2,\gamma,c)=\frac{2(\omega_1-\omega_2)}{c(x)^2}\imath\nonumber\\
&&\cdot\Bigl(
\nabla\hat{\phi}_1\cdot\overline{\nabla\hat{\phi}_2} 
+\omega_1\omega_2 \frac{\gamma(x)-1}{2c(x)^2}\hat{\phi}_1\overline{\hat{\phi}_2}\Bigr)\nonumber
\end{eqnarray}
which nicely illustrates the physical fact that the propagating wave described by $\psi$ is concentrated at the difference frequency $\omega_1-\omega_2$.

We mention in passing that in fact also in the harmonic ansatz for $\phi_1$, $\phi_2$ taking the real part would be demanded by physics. This would lead to certain (actually higher frequency) correction terms, that we neglect here, though, as they are not relevant for reconstructions.

\paragraph{Boundary conditions}
We consider a bounded computational domain $\Omega$, where the excitation surfaces $\Sigma_k$ are part of the boundary $\Sigma_k\subseteq\partial\Omega$ and the rest of $\partial\Omega$ is subject to impedance
boundary conditions in order to damp reflected waves
\[
\partial_\nu \hat{\phi}_k=-\imath\sigma_k \hat{\phi}_k  \mbox{ on }\partial\Omega\setminus\Sigma_k\,, 
\quad \partial_\nu \hat{\psi}=-\imath\sigma \hat{\psi} \mbox{ on }\partial\Omega\,.
\]
with nonnegative $L^\infty$ impedance coefficients $\sigma,\sigma_k$ that are bounded away from zero on an open subset of $\partial\Omega$ or $\partial\Omega\setminus\Sigma_k$, respectively.
\inarxiv{Note that with the choice  
$\sigma_k=\omega_k\sqrt{\kappa_0}$,
$\sigma=(\omega_1-\omega_2)\sqrt{\kappa_0}$,
these would be first order absorbing boundary conditions; however, later on in the definition of the operator $\mathcal{A}_c$ we wish to avoid explicit frequency dependence.
An alternative scenario that allows to work on a bounded domain $\Omega$ as well is to use a perfectly matched layer PML, (see. e.g., \cite[Section 5.5]{MKbook} and the references therein,) making the replacements 
$
\Delta\, \leftrightarrow \, \nabla\cdot(D(x)\nabla)$,
$\frac{\omega_k^2}{c(x)^2} \, \leftrightarrow \, d(x) \frac{\omega_k^2}{c(x)^2}
$
in the above Helmholtz equations on an augmented domain 
$\Omega=\Omega_{\rm acou}\cup\Omega_{\rm PML}$ with space dependent coefficients $D$ (matrix valued) and $d$; 
for details see, e.g., \cite{MKbook}.
Note that the real part of $D$  and $d$ is close to unity also in the PML region.
The boundary condition on the outer boundary can then be set to homogeneous Neumann.
}
\paragraph{Measurements}
The pressure data taken at the receiver array can, via the identity 
\[
\varrho\partial_t\psi=p
\]
be expressed by an observation operator 
\begin{equation}\label{eq:C}
C:(\hat{\phi}_1,\hat{\phi}_1,\hat{\psi})\mapsto \imath(\omega_1-\omega_2)\mbox{tr}_\Gamma\hat{\psi} 
\end{equation}
where $\Gamma$ is a manifold representing the receiver array and lying inside the acoustic domain $\Omega$.

\paragraph{Inverse problem}
The inverse problem of vibro-acoustography consists of determining the spatially varying coefficients $c$ and $\gamma$ from observations \eqref{eq:C} of the low frequency wave field.
We assume that $c$ is known on the outer boundary and needs to be reconstructed only in a subdomain (region of interest) $\widetilde{\Omega}\subseteq\Omega$ of the computational domain. With a slight abuse of notation we write
\[\frac{1}{c^2}=:\kappatil=\kappa_0+\chi_{\widetilde{\Omega}}\kappa\,, \quad \frac{\gamma-1}{2c^4}=:\gammatil=\gamma_0+\chi_{\widetilde{\Omega}}\gamma\]
where the background $\kappa_0$, $\gamma_0\in L^\infty(\widetilde{\Omega})$, $\kappa_0(x)\geq\underline{\kappa}>0$, and the subdomain $\tilde{\Omega}$ are known and $\chi_{\widetilde{\Omega}}$ is the extension by zero operator from $\widetilde{\Omega}$ to $\Omega$, defined by $(\chi_{\widetilde{\Omega}}\kappa)(x)=\kappa(x)$ for $x\in\widetilde{\Omega}$ and zero else.
Therewith our aim is to recover $\kappa,\gamma\in L^2(\widetilde{\Omega})$, in the weak form of \eqref{eq:phikhat0}, \eqref{eq:psihat0}
\begin{equation} \label{eq:A0}
\begin{aligned}
&0=\langle A(\hat{\phi}_1,\hat{\phi}_2,\hat{\psi},\kappa,\gamma),(v_1,v_2,w)\rangle:=\\
&\Re\Bigl((1-\imath)\Bigl(\sum_{k=1}^2\int_{\Omega} (-\omega_k^2\kappatil\hat{\phi}_k\overline{v}_k+\nabla\hat{\phi}_k\cdot\nabla \overline{v}_k)\, dx\\
&+\int_{\Omega} (-(\omega_1-\omega_2)^2\kappatil\hat{\psi}\overline{w}+\nabla\hat{\psi}\cdot\nabla \overline{w})\, dx\\
&+\imath\int_{\partial\Omega\setminus\Sigma}\sigma_k\hat{\phi}_k \overline{v}_k \, ds
-\int_\Sigma \hat{g}_k \overline{v}_k\, ds 
+\imath\int_{\partial\Omega}\sigma\hat{\psi} \overline{w} \, ds\\
&-2(\omega_1-\omega_2) \imath \\
&\quad\cdot\int_\Omega \left(
\kappatil\nabla\hat{\phi}_1\cdot\overline{\nabla\hat{\phi}_2} 
+\omega_1\omega_2 \gammatil \hat{\phi}_1\overline{\hat{\phi}_2}\right)\,\overline{w}\, dx \Bigr)\Bigr)
\\&\mbox{for all } v_1,v_2,w\in H^1(\Omega;\mathbb{C})
\end{aligned}
\end{equation}
(it suffices to take the real part here since $v_1,v_2,w$ vary over complex valued functions).
This is the weak form of 
\begin{eqnarray}
&-\omega_k^2\kappatil \hat{\phi}_k-\Delta \hat{\phi}_k=0 \mbox{ in }\Omega 
\nonumber\\
&\partial_\nu \hat{\phi}_k=\hat{g}_k \mbox{ on } \Sigma_k\,,\quad k\in\{1,2\}
\label{phihatk}\\
&-(\omega_1-\omega_2)^2\kappatil \hat{\psi}-\Delta \hat{\psi}=f(\kappatil,\gammatil,\hat{\phi}_1,\hat{\phi}_2)\mbox{ in }\Omega
\label{psihat}\\
& \quad f(\kappatil,\gammatil,\hat{\phi}_1,\hat{\phi}_2)=2\imath(\omega_1-\omega_2)
\nonumber\\
&\cdot\Bigl(\kappatil
\nabla\hat{\phi}_1\cdot\overline{\nabla\hat{\phi}_2} 
+\omega_1\omega_2 \gammatil\hat{\phi}_1\overline{\hat{\phi}_2}\Bigr)
\nonumber
\end{eqnarray}
with homogeneous impedance boundary conditions on (the rest of) $\partial\Omega$, which we also tacitly assume to hold in the following.

In Section~\ref{sec:forward} we will prove that for every $\kappa,\gamma\in L^2(\widetilde{\Omega})$, there exists a unique solution $\hat{\phi}_1,\hat{\phi}_2,\hat{\psi}$ of the operator equation $A(\hat{\phi}_1,\hat{\phi}_2,\hat{\psi},\kappa,\gamma)=0$ in appropriate function spaces, such that also the observation operator $C$ according to \eqref{eq:C} can be applied and yields an element of $L^2(\Gamma)$.

This justifies the use of the function spaces 
\begin{equation}\label{defXY}
X=L^2(\widetilde{\Omega})\times L^2(\widetilde{\Omega})\,, \quad  Y=L^2(\Gamma),
\end{equation}
to define the forward operator 
\begin{equation}\label{defF}
\begin{aligned}
&F:X\to Y, \quad \kappagamma\mapsto F(\kappagamma)=\imath(\omega_1-\omega_2)\mbox{tr}_\Gamma\hat{\psi}\\ &\mbox{ where $\phi_1,\phi_2,\psi$ solve \eqref{eq:A0} }
\end{aligned}
\end{equation}
and write the inverse problem in reduced form as 
\begin{equation}\label{Fxy}
F(\kappagamma)=y
\end{equation}
where $y\in L^2(\Gamma)$ is the pressure distribution measured at the receiver array. 
Concerning the choice of spaces \eqref{defXY}, working in $L^2$ spaces makes definition of methods most convenient. This is on one hand due to their Hilbert space structure, on the other hand due to the fact that no derivatives are involved, which avoids having to solve additional PDEs for evaluating the adjoint operator.

Alternatively, using the model and observation operators $A$ and $C$ defined in \eqref{eq:C}, \eqref{eq:A0}, we may write the inverse problem as an all-at-once system for the parameters $\kappagamma$ and the states $u:=(\hat{\phi}_1,\hat{\phi}_2,\hat{\psi})$ as 
\begin{equation}\label{AuxCu}
\begin{aligned}
&A(u,\kappagamma)=0\\
&Cu=y
\end{aligned}
\end{equation}

The two formulations are related via the identity $F=C\circ S$, where the parameter-to-state map $S:\kappagamma\to u:=(\hat{\phi}_1,\hat{\phi}_2,\hat{\psi})$ is implicitly defined by the identity
\begin{equation}\label{eq:S}
A(S(\kappagamma),\kappagamma)=0\,.
\end{equation}
\inarxiv{
Thus, using the forward operator $F$ requires an analysis of the operator $S$.
}

\section{Forward problem and function space setting}\label{sec:forward}

In the following, function spaces such as $L^2(\Omega;\mathbb{C})$ or $H^1(\Omega;\mathbb{C})$ will be regarded as spaces of functions with values in $\mathbb{C}$, but treated as real Hilbert spaces with a real valued inner product, e.g. $(v,w)_{L^2(\Omega)}=\Re(\int_\Omega v\overline{w}\, dx)$.
The $L^2$ space of real valued functions will simply denoted by $L^2(\Omega)$.

Consider the Laplace operator equipped with impedance boundary conditions, defined in its weak form by 
\[
\begin{aligned}
&\langle D_\sigma(\hat{\psi},w\rangle:=B_\sigma(\hat{\psi},w)\\
&:=
\Re\Bigl((1-\imath)\Bigl(\int_{\Omega} \nabla\hat{\psi}\cdot\nabla \overline{w})\, dx+\imath\int_{\partial\Omega}\sigma\hat{\psi} \overline{w} \, ds\Bigr)\Bigr), \\ 
&\forall w\in H^1(\Omega)\,.
\end{aligned}
\]
Here $B_\sigma$ is a symmetric, bounded and coercive bilinear form on on $H^1(\Omega;\mathbb{C})$ by the identity 
\[
B_\sigma(\hat{\psi},\hat{\psi})= 
\int_{\Omega} \vert\nabla\hat{\psi}\vert^2\, dx+\int_{\partial\Omega}\sigma\vert\hat{\psi}\vert^2 \, ds
\]
and Poincar\'{e}'s inequality. Thus, by the Lax-Milgram Lemma, $D_\sigma:H^1(\Omega;\mathbb{C})\to H^1(\Omega;\mathbb{C})^*$ is boundedly invertible and its inverse is compact as an operator from $L^2_{\kappatil}(\Omega;\mathbb{C})$ into itself, where $L^2_{\kappatil}(\Omega;\mathbb{C})$ is the weighted $L^2$ space with weight function $\kappatil\in L^2(\Omega)$, $\kappatil\geq0$ almost everywhere. Thus, by spectral theory for compact operators, $D_\sigma$ has a countable sequence of positive real eigenvalues tending to infinity, which we will denote by $\{\lambda^\sigma_n\, : n\in \mathbb{N}\}$.
Likewise, the eigenvalues of the operators $D_{\sigma_k}$ defined by the Laplacian on $\Omega$ with impedance boundary conditions (coefficient $\sigma_k$) on $\partial\Omega\setminus\Sigma_k$ are given by the countable set $\{\lambda^{\sigma_k}_n\, : n\in \mathbb{N}\}$, $k\in\{1,2\}$.
Thus \eqref{eq:A0} is uniquely solvable provided $\hat{g}_k\in H^{-\frac12}(\Sigma_k)$ and $\omega_k\notin \{\lambda^{\sigma_k}_n\, : n\in \mathbb{N}\}$, $k\in\{1,2\}$, $\omega_1-\omega_2\notin\{\lambda^\sigma_n\, : n\in \mathbb{N}\}$.

Higher regularity (actually only higher summability) can be achieved
under the additional assumption $\hat{g}_k\in (W^{1-\frac{1}{q},\frac{q}{q-1}}(\partial\Omega\setminus\Sigma_k))^*\subseteq(\mbox{tr}_{\Sigma_k}(W^{1,\frac{q}{q-1}}(\Omega)))^*$, meaning that the linear map $v\mapsto \int_\Sigma\hat{g}_k v\, ds$ lies in $(W^{1,\frac{q}{q-1}}(\Omega)))^*$.
Therefore according to elliptic regularity (e.g. \cite[Theorem 7.7]{Troeltzsch:2010}), \eqref{phihatk}
admits weak solutions $\hat{\phi}_k\in W^{1,q}(\Omega;\mathbb{C})$, $k\in\{1,2\}$.
Thus, the right hand side of \eqref{psihat} has the following regularity.
From $\kappatil\in L^2(\Omega)$ and $\nabla\hat{\phi}_1$, $\overline{\nabla\hat{\phi}_2} \in L^{q}(\Omega)$ we conclude by H\"older's inequality 
\begin{equation}\label{Hoelderabc}
\begin{aligned}
&\|a\,b\,\overline{c}\|_{L^r}\leq \|a\|_{L^2}\|b\|_{L^{\frac{4r}{2-r}}}\|c\|_{L^{\frac{4r}{2-r}}}\\
&\mbox{ for any }a\in L^2(\Omega), \ b,\,c\in L^{\frac{4r}{2-r}}(\Omega,\mathbb{C})
\end{aligned}
\end{equation}
that $\kappatil\nabla\hat{\phi}_1\cdot\overline{\nabla\hat{\phi}_2} \in L^r(\Omega)\subseteq W^{-1,p}(\Omega)$, provided
\begin{equation}\label{eq:qrp}
r\leq\min\left\{2,\frac{q}{2}\right\}\mbox{ and }
\frac{2r}{2-r}\leq\frac{q}{2} \mbox{ and }
1-\frac{d}{p^*}\geq -\frac{d}{r^*}
\end{equation}
where $p^*=\frac{p}{p-1}$ denotes the dual index.
This regularity (and even more) also holds true for the second quadratic term $\gammatil\hat{\phi}_1\overline{\hat{\phi}_2}$, with $\gammatil\in L^2(\Omega)$. Thus we conclude $\psi\in W^{1,p}(\Omega)$ (cf. \cite[Theorem 7.7]{Troeltzsch:2010}) and hence, by the Trace Theorem,  $\mbox{tr}_\Gamma\psi\in W^{1-\frac{1}{p},p}(\Gamma)\subseteq L^2(\Gamma)$ provided
\begin{equation}\label{eq:p}
1-\frac{1}{p}-\frac{d-1}{p}\geq -\frac{d-1}{2}\,.
\end{equation}
It is readily checked that conditions \eqref{eq:qrp}, \eqref{eq:p} can be satisfied for $\Omega\subseteq\mathbb{R}^d$, 
by choosing, 
\[
\begin{aligned}
&\frac{2d}{d+2}\leq p^*\leq\frac{2d}{d-1}, \\
&2\leq r^*\leq\frac{d p^*}{d-p^*} (\mbox{or }p^*\geq d), \quad
q\geq \frac{4r}{2-r},
\end{aligned}
\]
that is, in the physically relevant case $d\leq 3$, e.g. 
$p=\frac32$, $r=1$, $q=4$.

Thus we have proven
\begin{theorem}\label{thm:forward}
Let $\hat{g}_k\in H^{-\frac12}(\Sigma_k)$ and $\omega_k\notin \{\lambda^{\sigma_k}_n\, : n\in \mathbb{N}\}$, $k\in\{1,2\}$, $\omega_1-\omega_2\notin\{\lambda^\sigma_n\, : n\in \mathbb{N}\}$ the sets of eigenvalues of the Laplacians $D_{\sigma_k}$, $D_\sigma$ with impedance boundary conditions.

Then the parameter-to-state map $S:\mathcal{D}(F)\to W^{1,q}(\Omega;\mathbb{C})^2\times W^{1,p}(\Omega;\mathbb{C})$, and the forward operator $F:\mathcal{D}(F)\to L^2(\Gamma)$ are well-defined by \eqref{eq:C}, \eqref{eq:A0}, \eqref{defF}, \eqref{eq:S} on 
$\mathcal{D}(F)=\{\kappa\in L^2(\widetilde{\Omega})\, : \, \kappatil\geq0\mbox{ a.e. }\}\times L^2(\widetilde{\Omega})$.
\end{theorem}

The domain $\mathcal{D}(F)$ has empty interior with respect to the $L^2$ topology and this prevents applicability of convergence results for the Newton and gradient methods to be discussed below.
To avoid this, we restrict $F$ to an open ball around a strictly positive $L^\infty$ function $\kappa_0>0$ (e.g., the background)
\begin{equation}\label{DF}
\begin{aligned}
\tilde{\mathcal{D}}(F)=&\{\kappa\in L^2(\widetilde{\Omega})\, : \, \|\kappatil-\kappa_0\|_{L^2}\leq\rho\}\\
&\times L^2(\widetilde{\Omega})
\end{aligned}
\end{equation}
for $\rho$ sufficiently small, and apply a fixed point argument to obtain well-definedness of $F$ on $\tilde{\mathcal{D}}(F)$, see, e.g., \cite{HNS95}.

For use in Newton and gradient type methods we also need differentiability of $F$. It sufficies to prove that the parameter-to-state map $S$ is differentiable, since $F=C\circ S$ with $C$ being a bounded linear operator.
It is straightforward to see that for
$(\hat{\phi}_1,\hat{\phi}_2,\hat{\psi}):=S(\kappagamma)$, 
$(\hat{\phi}_1^+,\hat{\phi}_2^+,\hat{\psi}^+):=S(\kappagamma+\delta\kappagamma)$
the difference $(d\hat{\phi}_1,d\hat{\phi}_2,d\hat{\psi}):=S(\kappagamma+\delta\kappagamma)-S(\kappagamma)$ satisfies the weak form of 
\begin{equation}\label{phihatkpsihat_diff}
\begin{aligned}
&-\omega_k^2\kappatil \, d\hat{\phi}_k-\Delta \, d\hat{\phi}_k=\omega_k^2\,\underline{\widetilde{\delta\kappa}}\, \hat{\phi}_k^+ \mbox{ in }\Omega \\ 
&\partial_\nu \, d\hat{\phi}_k=0 \mbox{ on } \Sigma_k\,,\quad k\in\{1,2\}
\\
&-(\omega_1-\omega_2)^2\kappatil \, d\hat{\psi}-\Delta \, d\hat{\psi}=f_d
\mbox{ in }\Omega\\
& \quad f_d= (\omega_1-\omega_2)^2\kappatil \, \hat{\psi}^+ +2\imath(\omega_1-\omega_2)\\
&\cdot\Bigl(\underline{\widetilde{\delta\kappa}}\,\nabla\hat{\phi}_1^+\cdot\overline{\nabla\hat{\phi}_2^+} 
+\kappatil(\nabla\, d\hat{\phi}_1\cdot\overline{\nabla\hat{\phi}_2^+} 
+\nabla\hat{\phi}_1\cdot\overline{\nabla\, d\hat{\phi}_2})\\
&\qquad + \omega_1\omega_2 \left(
\underline{\widetilde{\delta\gamma}}\,\hat{\phi}_1^+\overline{\hat{\phi}_2^+}
+\gammatil( d\hat{\phi}_1\overline{\hat{\phi}_2^+}
+\hat{\phi}_1\,\overline{d\hat{\phi}_2})
\right)\Bigr)
\end{aligned}
\end{equation}
with 
$\underline{\widetilde{\delta\kappa}}=\chi_{\widetilde{\Omega}}\delta\kappa$,
$\underline{\widetilde{\delta\gamma}}=\chi_{\widetilde{\Omega}}\delta\gamma$,
and therefore, formally
$(\delta\hat{\phi}_1,\delta\hat{\phi}_2,\delta\hat{\psi}):=S'(\kappagamma)\delta\kappagamma$
solves
\begin{equation}\label{phihatkpsihat_deriv}
\begin{aligned}
&-\omega_k^2\kappatil \, \delta\hat{\phi}_k-\Delta \, \delta\hat{\phi}_k=\omega_k^2\,\underline{\widetilde{\delta\kappa}}\, \hat{\phi}_k \mbox{ in }\Omega \\
&\partial_\nu \, \delta\hat{\phi}_k=0 \mbox{ on } \Sigma_k\,,\quad k\in\{1,2\}
\\
&-(\omega_1-\omega_2)^2\kappatil \, \delta\hat{\psi}-\Delta \, \delta\hat{\psi}=
f_\delta\mbox{ in }\Omega\\
& \quad f_\delta=(\omega_1-\omega_2)^2\kappatil \, \hat{\psi} +2\imath(\omega_1-\omega_2)\\
&\cdot\Bigl(\underline{\widetilde{\delta\kappa}}\,\nabla\hat{\phi}_1\cdot\overline{\nabla\hat{\phi}_2} 
+\kappatil(\nabla\, \delta\hat{\phi}_1\cdot\overline{\nabla\hat{\phi}_2} 
	+\nabla\hat{\phi}_1\cdot\overline{\nabla\, \delta\hat{\phi}_2})\\
&\qquad + \omega_1\omega_2 \left(
\underline{\widetilde{\delta\gamma}}\,\hat{\phi}_1\overline{\hat{\phi}_2}
+\gammatil( \delta\hat{\phi}_1\overline{\hat{\phi}_2}
	+\hat{\phi}_1\,\overline{\delta\hat{\phi}_2})
\right)\Bigr)\,.
\end{aligned}
\end{equation}
Hence the first order Taylor remainder 
$(\phiTay_1,\phiTay_2,\psiTay):=S(\kappagamma+\delta\kappagamma)-S(\kappagamma)-S'(\kappagamma)\delta\kappagamma$
obeys 
\begin{equation}\label{phihatkpsihat_remainder}
\begin{aligned}
&-\omega_k^2\kappatil \, \phiTay_k -\Delta \, \phiTay_k=\omega_k^2\,\underline{\widetilde{\delta\kappa}}\, d\hat{\phi}_k \mbox{ in }\Omega \\
&\partial_\nu \, \phiTay_k=0 \mbox{ on } \Sigma_k\,,\quad k\in\{1,2\}\\
&-(\omega_1-\omega_2)^2\kappatil \, \psiTay-\Delta \, \psiTay=f_{rest}\mbox{ in }\Omega\\
& \quad f_ {rest}=
(\omega_1-\omega_2)^2\kappatil \, d\hat{\psi} +2\imath(\omega_1-\omega_2)\\
&\cdot\Bigl(\underline{\widetilde{\delta\kappa}}\,
(\nabla\, d\hat{\phi}_1\cdot\overline{\nabla\hat{\phi}_2^+} 
+\nabla\hat{\phi}_1\cdot\, \overline{d\nabla\hat{\phi}_2})\\
&\ +\kappatil(\nabla\, \phiTay_1\cdot\overline{\nabla\hat{\phi}_2}
+\nabla\, d\hat{\phi}_1\cdot\overline{\nabla\, d\hat{\phi}_2} 
	+\nabla\hat{\phi}_1\cdot\overline{\nabla \phiTay_2})\\
&\ + \omega_1\omega_2 \Bigl(
\underline{\widetilde{\delta\gamma}}\,
(d\hat{\phi}_1\overline{\hat{\phi}_2^+} 
+\hat{\phi}_1\, \overline{d\hat{\phi}_2})\\
&\qquad\qquad+\gammatil(\phiTay_1\overline{\hat{\phi}_2}
+\,d\hat{\phi}_1\,d\overline{\hat{\phi}_2}
	+\hat{\phi}_1\overline{\phiTay_2})
\Bigr)\Bigr).
\end{aligned}
\end{equation}
Here we have used the identities
\[
\begin{aligned}
&(a+\, \delta a) (b+db) (c+\,dc) - abc\\
&\quad= \delta a\, (b+db)(c+\,dc)\, + \, a\, db \, (c+dc)\, +\, a\,b\,dc\\
&(a+\, \delta a) (b+db) (c+\,dc) - abc\\ 
&\quad-(\delta a\, b\, c\, + \, a\, \delta b\, c\, +\, a\,b\,\delta c)\\
&\quad=\delta a\,(db\, (c+dc)\, + \, b\, dc) \\
&\qquad+ a[(db-\delta b)c+db\,dc+b(dc-\delta c)].
\end{aligned}
\]
Regularity arguments as in the proof of Theorem~\ref{thm:forward} lead to estimates of the form 
\[
\begin{aligned}
&\|d\hat{\phi}_k\|_{W^{1,q}}\leq C \|\underline{\widetilde{\delta\kappa}}\, \hat{\phi}_k^+\|_{(W^{1,q^*})^*}\\
&\|d\hat{\psi}\|_{W^{1,p}}\leq C \|f_d\|_{(W^{1,p^*})^*}\\
&\|\delta\hat{\phi}_k\|_{W^{1,q}}\leq C \|\underline{\widetilde{\delta\kappa}}\, \hat{\phi}_k\|_{(W^{1,q^*})^*}\\
&\|\delta\hat{\psi}\|_{W^{1,p}}\leq C \|f_\delta\|_{(W^{1,p^*})^*}\\
&\|\phiTay_k\|_{W^{1,q}}\leq C \|\underline{\widetilde{\delta\kappa}}\, d\hat{\phi}_k\|_{(W^{1,q^*})^*}\\
&\|\psiTay\|_{W^{1,p}}\leq C \|f_{rest}\|_{(W^{1,p^*})^*}\,,
\end{aligned}
\]
where $f_d$, $f_\delta$, $f_{rest}$ can be estimated by the same H\"older inequalities and Sobolev embeddings as those used for the proof of Theorem \ref{thm:forward}. 

This proves Fr\'{e}chet differentiability.
\begin{theorem}\label{thm:forward_diff}
Under assumptions of Theorem~\ref{thm:forward}, the parameter-to-state map $S$ and the forward operator $F$ are Fr\'{e}chet differentiable on $\tilde{\mathcal{D}}(F)$ as defined in \eqref{DF} with respect to the $L^2$ topology in preimage space, as mappings to $W^{1,q}(\Omega;\mathbb{C})^2\times W^{1,p}(\Omega;\mathbb{C})$ and $L^2(\Gamma)$, respectively.
\end{theorem}

Concerning further convergence conditions for Newton and gradient type methods, cf. e.g. \cite{KNSbook:2008}, we briefly comment on the tangential cone condition
\begin{equation}\label{tcc}
\begin{aligned}
&\|F(\kappagamma+\delta\kappagamma)-F(\kappagamma)-F'(\kappagamma)\delta\kappagamma\|\\
&\quad\leq c_{tc}\|F(\kappagamma+\delta\kappagamma)-F(\kappagamma)\|.
\end{aligned}
\end{equation}
In case of known speed of sound $c$, when we seek to identify $\gamma=\gamma(x)$ only, 
the inverse problem becomes an inverse source problem, see \eqref{eq:psihat_hf} below, and is therefore affinely linear, thus trivially satisfying \eqref{tcc} with $c_{tc}=0$.
Conversely, if $c=c(x)$ is to be determined, the inverse probems is closely related to the well-known and well-investigated model problem of recovering the potential $c$ in the Schr\"odinger equation $-\Delta u + c u=0$. This is known to satisfy the tangential cone condition only in case of complete observations of $u$ on all of $\Omega$ \cite{HNS95}. Thus \eqref{tcc} cannot be expected to be verifiable in our boundary observation setting.

In the definition of gradient type methods (and also in the implementation of Newton type methods) we will need the adjoint of $F'(\kappagamma)$, which we therefore derive here.
First of all, note that by $F=C\circ S$ with $S$ defined by \eqref{eq:S} and the Implicit Function Theorem we can write 
\begin{equation}\label{eq:SpFpKL}
S'(\kappagamma)=-K^{-1} L\,, \quad 
F'(\kappagamma)=-CK^{-1} L\,, 
\end{equation}
where  
\begin{equation}\label{eq:KL}
K=\frac{\partial A}{\partial u}(u^{(n)},\kappagamma^{(n)})\,, \quad 
L=\frac{\partial A}{\partial \kappagamma}(u^{(n)},\kappagamma^{(n)})
\end{equation}
are the linearizations of the operator $A$ from \eqref{eq:A0} with respect to the states and the parameters, respectively. They are given by
\begin{equation}\label{eq:Kmod}
\begin{aligned}
&\langle K\, (\delta\hat{\phi}_1,\delta\hat{\phi}_2,\delta\hat{\psi}),(v_1,v_2,w)\rangle=\\
&\Re\Bigl((1-\imath)\Bigl(\sum_{k=1}^2\int_{\Omega} (-\omega_k^2\kappatil\, \delta\hat{\phi}_k\overline{v}_k+\nabla\ \delta\hat{\phi}_k\cdot\nabla \overline{v}_k)\, dx\\
&+\int_{\Omega} (-(\omega_1-\omega_2)^2\kappatil\, \delta\hat{\psi}\overline{w}+\nabla\, \delta\hat{\psi}\cdot\nabla \overline{w})\, dx\\
&+\imath\int_{\partial\Omega\setminus\Sigma}\sigma_k\, \delta\hat{\phi}_k \overline{v}_k \, ds
+\imath\int_{\partial\Omega}\sigma\, \delta\hat{\psi} \overline{w} \, ds\\
&-\int_\Omega 2(\omega_1-\omega_2) \imath \Bigl(
\kappatil
\Bigl(\nabla\, \delta\hat{\phi}_1\cdot\overline{\nabla\hat{\phi}_2}
+\nabla\hat{\phi}_1\cdot\overline{\nabla\, \delta\hat{\phi}_2}\Bigr)\\
&\quad+\omega_1\omega_2 \gammatil \Bigl(
\, \delta\hat{\phi}_1\,\overline{\hat{\phi}_2}
+\hat{\phi}_1\overline{\, \delta\hat{\phi}_2}
\Bigr)
\Bigr)\Bigr)\,\overline{w}\, dx \Bigr)
\end{aligned}
\end{equation}
\begin{equation}\label{eq:Lmod}
\begin{aligned}
&\langle L\, (\delta\kappa,\delta\gamma),(v_1,v_2,w)\rangle=\\
&\Re\Bigl((1-\imath)\Bigl(\sum_{k=1}^2\int_{\Omega} -\omega_k^2 \underline{\widetilde{\delta\kappa}} \,\hat{\phi}_k\overline{v}_k\, dx\\
& +\int_{\Omega} -(\omega_1-\omega_2)^2 \underline{\widetilde{\delta\kappa}} \,\hat{\psi}\overline{w}\, dx\\
&-2(\omega_1-\omega_2) \imath \\
&\cdot\int_\Omega \left(
\underline{\widetilde{\delta\kappa}} \,\nabla\hat{\phi}_1\cdot\overline{\nabla\hat{\phi}_2} 
+\omega_1\omega_2 \, \underline{\widetilde{\delta\gamma}}\, \hat{\phi}_1\overline{\hat{\phi}_2}\right)\,\overline{w}\, dx \Bigr)\Bigr)
\end{aligned}
\end{equation}
for any $v_1,v_2,w\in H^1(\Omega;\mathbb{C})$.
The identity \eqref{eq:SpFpKL} with \eqref{eq:KL} can also be used to determine the adjoint operator $F'(\kappagamma)^*=-(CK^{-1}L)^*$ as a Hilbert space adjoint in $L^2$. 
To this end, for given $r\in L^2(\Gamma)$ we want to find $\xizeta:=(CK^{-1}L)^*r$ such that
\[
\begin{aligned}
&\langle CK^{-1}L\,\delta\kappagamma,r\rangle_{L^2(\Gamma)} = \langle \delta\kappa,\xi\rangle_{L^2(\widetilde{\Omega})}
+\langle \delta\gamma,\zeta\rangle_{L^2(\Omega)} \\
&\mbox{ for all }\delta\kappagamma\in L^2(\widetilde{\Omega})\times L^2(\Omega)\,.
\end{aligned}
\]
We introduce the auxiliary variables $(\widetilde{\phi}_1,\widetilde{\phi}_2,\widetilde{\psi}):=K^{-1}L\,\delta\kappagamma$, which allows us to use the identity 
\begin{equation}\label{eq:stateeq}
\begin{aligned}
&\langle K(\widetilde{\phi}_1,\widetilde{\phi}_2,\widetilde{\psi}),(v_1,v_2,w)\rangle 
=\langle L\delta\kappagamma,(v_1,v_2,w)\rangle\\
&\mbox{ for all }(v_1,v_2,w)
\end{aligned}
\end{equation}
and define $(p_1,p_2,q)$ as the solution to the adjoint equation
\begin{equation}\label{eq:adjeq}
\begin{aligned}
&\langle K(\delta\hat{\phi}_1,\delta\hat{\phi}_2,\delta\hat{\psi}),(p_1,p_2,q)\rangle\\
&\qquad=\langle C(\delta\hat{\phi}_1,\delta\hat{\phi}_2,\delta\hat{\psi}),r\rangle_{L^2(\Gamma)} 
\\&
\mbox{ for all } (\delta\hat{\phi}_1,\delta\hat{\phi}_2,\delta\hat{\psi})\,.
\end{aligned}
\end{equation}
Using \eqref{eq:stateeq}, and \eqref{eq:adjeq} together with \eqref{eq:Kmod}, \eqref{eq:Lmod}, we get
\[
\begin{aligned}
&\langle C(\widetilde{\phi}_1,\widetilde{\phi}_2,\widetilde{\psi}),r\rangle_{L^2(\Gamma)} 
=\langle L\,\delta\kappagamma,(p_1,p_2,q)\rangle\\
&=-\Re\Bigl(\sum_{k=1}^2\int_{\widetilde{\Omega}} \omega_k^2 \,\delta\kappa \,\hat{\phi}_k^{(n)}\overline{p}_k\, dx\\ 
&\qquad+\int_{\widetilde{\Omega}} (\omega_1-\omega_2)^2 \,\delta\kappa \,\hat{\psi}^{(n)}\overline{q}\, dx\\
&\qquad+2(\omega_1-\omega_2) \imath \\
&\hspace*{-0.5cm}\cdot\int_{\widetilde{\Omega}} \left(
\delta\kappa \,\nabla\hat{\phi}_1^{(n)}\cdot\overline{\nabla\hat{\phi}_2^{(n)}} 
+\omega_1\omega_2 \, \delta\gamma\, \hat{\phi}_1^{(n)}\overline{\hat{\phi}_2^{(n)}}\right)\,\overline{q}\, dx \Bigr)\\
&=\langle \delta\kappa,\xi\rangle_{L^2(\widetilde{\Omega})}
+\langle \delta\gamma,\zeta\rangle_{L^2(\widetilde{\Omega})}
\end{aligned}
\]
for 
\[
\begin{aligned}
&\xi=-\Re\Bigl(\sum_{k=1}^2 \omega_k^2 \hat{\phi}_k^{(n)}\overline{p}_k\\
&+(\omega_1-\omega_2)^2 \hat{\psi}^{(n)}\overline{q}
+2(\omega_1-\omega_2) \imath \nabla\hat{\phi}_1^{(n)}\cdot\overline{\nabla\hat{\phi}_2^{(n)}} \overline{q} \Bigr)\vert_{\widetilde{\Omega}}
\\
&\zeta=-\Re\Bigl(2(\omega_1-\omega_2) \imath \omega_1\omega_2 \hat{\phi}_1^{(n)}\overline{\hat{\phi}_2^{(n)}}\overline{q}\Bigr)\vert_{\widetilde{\Omega}}.
\end{aligned}
\]
Thus we end up with an explicit expression for $\xizeta:=F'(\kappagamma)^*r$.
For this purpose, the adjoint states have to be computed as solutions to \eqref{eq:adjeq}, that is, the weak form of 
\[
\begin{aligned}
&-(\omega_1-\omega_2)^2\kappatil q-\Delta q=0 \mbox{ in }\Omega\setminus\Gamma\,, \\ 
&\partial_\nu q=-\imath\sigma q \mbox{ on } \partial\Omega\,, \quad
\left[\partial_\nu q\right]=-\imath(\omega_1-\omega_2) r \mbox{ on }\Gamma 
\end{aligned}
\]
where $\left[\partial_\nu q\right]$ denotes the jump of the normal derivative over the interface $\Gamma$, as well as
\[
\begin{aligned}
&-\omega_k^2\kappatil p_k-\Delta p_k
=f_k \mbox{ in }\Omega\,, \\ 
&\partial_\nu p_k=-\imath\sigma_k p_k \mbox{ on } \partial\Omega\setminus\Sigma_k\,, \quad
\partial_\nu p_k=0 \mbox{ on } \Sigma_k\\
&f_k=-2(\omega_1-\omega_2)\imath\Bigl(
\nabla \bigl( \kappatil q\nabla \phi_{\not k}\bigr)
+\omega_1\omega_2 \gammatil q\phi_{\not{k}}\Bigr)
\end{aligned}
\]
with $\not{\!k}=\begin{cases}2\mbox{ for }k=1\\1\mbox{ for }k=2\end{cases}$.

\section{Uniqueness 
}\label{sec:uniqueness}
In case the speed of sound $c$ is known, reconstruction of $\gamma=\gamma(x)$ in the time domain \eqref{eq:phik}, \eqref{eq:psi} or the frequency domain formulation \eqref{eq:phikhat0}, \eqref{eq:psihat0} amounts to an inverse source problem.
Indeed, setting 
\[
\begin{aligned}
h(x,\omega_d;\omega_2)&=2\imath\omega_d\omega_2 \nabla\hat{\phi}(x,\omega_2+\omega_d)\cdot\overline{\nabla\hat{\phi}(x,\omega_2)}\\
\multip(x,\omega_d;\omega_2)&=2\imath 
(\omega_2+\omega_d)\omega_2 
\hat{\phi}(x,\omega_2+\omega_d)\overline{\hat{\phi}(x,\omega_2)}\\
\gamma'(x)&=\frac{\gamma(x)-1}{2c(x)^2}
\end{aligned} 
\]
and multiplying with $c^2$, we can write \eqref{eq:psihat0} as 
\begin{equation}\label{eq:psihat_hf}
-\omega_d^2 \hat{\psi}+\mathcal{A}_c \hat{\psi}= h(\omega_d;\omega_2)+ \multip(\omega_d;\omega_2) \gamma'
\end{equation}
Here we denote the difference frequency $\omega_1-\omega_2$ by $\omega_d$ and the solution of \eqref{eq:phikhat0} with boundary excitation $\hat{g}_k=\hat{g}(\omega)$ by $\hat{\phi}_k(\omega)$.
Note that the functions $\multip$ and $h$ are known from the known excitations $\hat{g}$.
Moroever, we denote by $\mathcal{A}_c$ 
\inarxiv{either of} 
the elliptic differential operator 
$-c^2\Delta$ with homogeneous impedance boundary conditions; 
\inarxiv{or $-\frac{c^2}{d}\nabla\cdot(D\nabla)$ with homogeneous Neumann boundary conditions. In both cases }
$\mathcal{A}_c$ is a selfadjoint nonnegative definite operator with respect to the weighted $L^2$ inner product with weight function $w=\frac{1}{c^2}$
\inarxiv{
or $w=\frac{d}{c^2}$, respectively}. 
By $\{(\lambda_k,(\varphi_j^k)_{j\in \indset_k}\, : \, k\in\mathbb{N}\}$ we denote the corresponding eigensystem, where in case of multiple eigenvalues we collect the eigenfunctions corresponding to $\lambda_k$ in the set $\{\varphi_j^k\, : \, j\in \indset_k\}$ with some finite index set $\indset_k$. (Note that in one space dimension, the eigenfunctions are single and so $\indset_k=\{1\}$.)
The requirements on $c$ for this purpose are 
\begin{equation}\label{eq:cond-c}
c,\tfrac{1}{c}\in L^\infty(\Omega)\,.
\end{equation}
Since the eigenfunctions form an orthonormal basis of $L^2_{w}$, we can expand $\hat{\psi}$ with respect to this basis 
\[
\hat{\psi}(x,\omega_d;\omega_2)
=\sum_{k=1}^\infty\sum_{j\in \indset_k}\langle \hat{\psi}(\omega_d;\omega_2),\varphi_j^k
\rangle_{L^2_{w}} \varphi_j^k(x)\,.
\]
This allows us to express the observations according to \eqref{eq:C} by 
\begin{equation}\label{eq:obs_eig}
\begin{aligned}
&y(x_0,\omega_d;\omega_2)=\hat{\psi}(x_0,\omega_d;\omega_2) \\
&=\sum_{k=1}^\infty\sum_{j\in \indset_k}\langle \hat{\psi}(\omega_d;\omega_2),\varphi_j^k
\rangle_{L^2_{w}} \varphi_j^k(x_0) \\ 
&x_0\in\Gamma, \quad \omega_d\in \freqset\,,
\end{aligned}
\end{equation}
where we assume that we can take observations for all difference frequencies in some set $\freqset $, while $\omega_2$ is fixed.
On the other hand, taking inner products of \eqref{eq:psihat_hf} with $\varphi_j^k$ and using the eigenvalue equation $\mathcal{A}_c\varphi_j^k=\lambda_j\varphi_j^k$ we obtain the identity
\begin{equation*}
\begin{aligned}
&\langle \hat{\psi}(\omega_d;\omega_2),\varphi_j^k\rangle_{L^2_{w}}\\
&= \frac{1}{-\omega_d^2+\lambda_j}
\langle h(\omega_d;\omega_2)+ \multip(\omega_d;\omega_2) \gamma',\varphi_j^k\rangle_{L^2_{w}}.
\end{aligned}
\end{equation*}
Combining this with \eqref{eq:obs_eig}  we get 
\begin{equation}\label{eq:obs_eig_sumj0}
\begin{aligned}
&\hspace*{-0.5cm}\tilde{y}(x_0,\omega_d;\omega_2)
=\sum_{k=1}^\infty\frac{1}{-\omega_d^2+\lambda_k}\sum_{j\in \indset_k}
\langle \multip(\omega_d;\omega_2) \gamma',\varphi_j^k\rangle_{L^2_{w}}\\
&x_0\in\Gamma, \quad \omega_d\in \freqset \,,
\end{aligned}
\end{equation}
where $\tilde{y}$ is the modified observation function
\[
\begin{aligned}
&\tilde{y}(x_0,\omega_d;\omega_2)\\
&=y(x_0,\omega_d;\omega_2)- (-\omega_d^2+\mathcal{A}_c)^{-1} h(\omega_d;\omega_2)\,,
\end{aligned}
\]
thus a known quantity.
In order to obtain from this the desired information on $\gamma'$, we assume that $\hat{g}(\omega)$ has been chosen such that $\multip$ factorizes into a frequency dependent and a space dependent part
\begin{equation}\label{eq:fab}
\multip(x,\omega_d;\omega_2)=a(\omega_d;\omega_2) b(x)
\end{equation}
so that \eqref{eq:obs_eig_sumj0} becomes 
\begin{equation}\label{eq:obs_eig_sumj}
\begin{aligned}
&\tilde{y}(x_0,\omega_d;\omega_2)
=\sum_{k=1}^\infty\frac{a(\omega_d;\omega_2)}{-\omega_d^2+\lambda_k}\sum_{j\in \indset_k}
\langle b \gamma',\varphi_j^k\rangle_{L^2_{w}}
\\& x_0\in\Gamma, \quad \omega_d\in \freqset \,.
\end{aligned}
\end{equation}
Both sides of this equality have sigularities at $\omega_d=\pm\sqrt{\lambda_\ell}$.
Thus, these poles provide the location of the eigenvalues of $\mathcal{A}_c$ and therewith some information on $c$ (see Remark~\ref{req:uniqueness_c} below).
Moreover, multiplying with $(\omega_d-\sqrt{\lambda_\ell})$ and taking the limit $\omega_d\to\sqrt{\lambda_\ell}$, we can extract the contribution due to the $\ell$th eigenfunction
\begin{equation}\label{eq:obs_eig_ell}
\begin{aligned}
&\lim_{\omega_d\to\sqrt{\lambda_\ell}} (\omega_d-\sqrt{\lambda_\ell})
\tilde{y}(x_0,\omega_d;\omega_2)\\
&= -\frac{a(\sqrt{\lambda_\ell};\omega_2)}{2\sqrt{\lambda_\ell}}\sum_{j\in \indset_k}
\langle b \gamma',\varphi_j^\ell\rangle_{L^2_{w}}
\quad x_0\in\Gamma
\end{aligned}
\end{equation}
For this to work out, we need to assume that
\begin{equation}\label{eq:cond_I}
\sqrt{\lambda_\ell}\mbox{ is an interior point of }\freqset \mbox{ for all }\ell\in\mathbb{N}. 
\end{equation}
Finally \eqref{eq:obs_eig_ell} allows to uniquely determine the coefficients $\langle b \gamma',\varphi_j^\ell\rangle_{L^2_{w}}$ in
\begin{equation}\label{gamma_recon}
\gamma'(x)=\frac{1}{b(x)}
\sum_{\ell=1}^\infty\sum_{j\in \indset_\ell} \langle b \gamma',\varphi_j^\ell\rangle_{L^2_{w}}
\end{equation}
provided $b$ vanishes nowhere and 
\begin{equation}\label{eq:cond_linindep}
\{\varphi_j^\ell\vert_\Gamma\,:\, j\in \indset_\ell\} \mbox{ is linearly independent} 
\end{equation}

Thus we have proven the following uniqueness result on recovery of $\gamma(x)$.
\begin{theorem}\label{thm:uniqueness}
Assume that $c$ is known and satisfies \eqref{eq:cond-c},
that $\freqset $ and $\Gamma$ are chosen such that \eqref{eq:cond_I}, \eqref{eq:cond_linindep}
hold, and that $\hat{g}(\omega_2+\omega_d)$, is chosen such that \eqref{eq:fab}
holds for all $\omega_d\in \freqset $ with $b\in L^\infty$.

Then $\gamma\in L^2(\Omega)$ is uniquely determined on the set $\{x\in\widetilde{\Omega}\, : \, b(x)\not=0\}$ by the observations
$y(x_0,\omega_d;\omega_2)=\hat{\psi}(x_0,\omega_d;\omega_2)$, $x_0\in\Gamma$, $\omega_d\in \freqset $.
\end{theorem}

\begin{remark}\label{rem:thm_uniqueness}
Obviously, if \eqref{eq:cond_I}, \eqref{eq:cond_linindep} only hold with $\mathbb{N}$ replaced by $\{1,\ldots,N\}$, we can recover the first $N$ coefficients of $b\gamma'$.

Note that no regularity assumptions with respect to $\omega_d$ need to be imposed here.

Condition \eqref{eq:cond_linindep} has been discussed in detail in \cite[Remark 4.1]{nonlinearity_imaging_fracWest}. It is trivially satisfied with $\Gamma$ containing a single point $\{x_0\}$ in one space dimension, since the eigenvalues of $\mathcal{A}_c$ are single then.
Moroever, it can be extended to higher space dimensions and geometric settings in which the eigenfunctions allow for separation of variables. A simple 2-d example is a disc with radius $r$, where using polar coordinates, the eigenfunctions can be written in terms of Bessel functions. A circle with almost any radius $r_*\in(0,r]$ can then be used as observation manifold $\Gamma$, as shown in  \cite[Remark 4.1]{nonlinearity_imaging_fracWest}.

To achieve the separability \eqref{eq:fab} of $\multip$ we supplement the boundary excitation $\hat{g}_k(\omega)$ by an interior one $f_g(\omega)$, which we view as an approximation of a source $\tilde{g}(\omega)\,\delta_\Sigma$ concentrated on $\Sigma$, cf., e.g., \cite{periodicWestervelt}.
The resulting equation for $\hat{\phi}_k(\omega)$
\[
\begin{aligned}
&-\omega^2\hat{\phi}_k -\frac{c^2}{d}\Delta \hat{\phi}_k = f_g(\omega)\mbox{ in }\Omega\\
&\partial_\nu \hat{\phi}_k=-\sigma_k\hat{\phi}_k\mbox{ on }\partial\Omega\setminus\Sigma_k\,,\
\partial_\nu \hat{\phi}_k=\hat{g}_k(\Omega)\mbox{ on }\Sigma_k
\end{aligned}
\]
then has a solution of the form $\hat{\phi}_k(x,\omega)=\tilde{a}(\omega)\tilde{b}(x)$
if, e.g., we choose $\tilde{b}$ such that $\Delta \tilde{b}=0$ in $\Omega$, $\partial_\nu \tilde{b}=-\sigma_k\tilde{b}$ on $\partial\Omega\setminus\Sigma$, and set 
$\hat{g}_k(\omega):=\tilde{a}(\omega)\partial_\nu \tilde{b}\vert_\Sigma$, 
$f_g(\omega):=-\omega^2\tilde{a}(\omega)\tilde{b}$.
\end{remark}

\begin{remark}
In case of constant sound speed $c$, uniqueness for the above inverse source problem for $\gamma(x)$ in the time domain formulation \eqref{eq:phik}, \eqref{eq:psi} from boundary observations under a space-time separability assumption (similar to the space-frequency one \eqref{eq:fab}) can be conclused from \cite[Theorem 7.4.2]{IsakovInvSourceBook}, provided $\Gamma$, $c$ and $T$ satisfy \cite[condition (1.2.11)]{IsakovInvSourceBook}, which is basically a condition on sufficient size of $\Gamma$ and $T$, depending on the speed $c$ of sound propagation.

Other related uniqueness results for $\gamma$ have been found recently in the context of nonlinearity imaging in \cite{nonlinearity_imaging_Westervelt,YamamotoBK:2021}.
\end{remark}

\begin{remark}\label{req:uniqueness_c}
Note that from the poles on both sides of \eqref{eq:obs_eig_sumj} we also obtain the eigenvalues of $\mathcal{A}_c$. According to Sturm-Liouville theory, applied as in \cite[Section 5.3]{fracPAT}, this uniquely determines $c(x)$ in one space dimension, provided we can take measurements at two different impedance values $\sigma$, $\tilde{\sigma}$. Note however, that we need the eigenfunctions of $\mathcal{A}_c$ for reconstructing $\gamma'(x)$ according to \eqref{gamma_recon}, so this only gives a uniqueness result for $c$ alone and no simultaneous uniqueness of $c$ and $\gamma$. Also, its restriction to the 1-d setting limits applicability to our experimental setting.

For uniqueness of $c=c(x)$ in higher space dimensions from boundary measurements, 
results on uniqueness of the space-dependent index of refraction $n(x)=\frac{c_0^2}{c(x)^2}$ in inverse scattering, e.g.,  \cite[Chapter 6]{KirschBuch} or of the potential in the Schr\"odinger equation \cite[Chapter 5]{IsakovPDE} are relevant.
Note however, that $c$ appears not only in the equation for the observed quantity $\psi$ but also governs the two excitation wave fields $\phi_1$, $\phi_2$ that enter the $\psi$ equation through a source term. This makes the uniqueness question for $c$ more involved than in the mentioned references.
\end{remark}

\begin{remark}
A proof of unique recovery of both $c$ and $\gamma$ is widely open and subject of future research. We will nevertheless in the remainder of this paper discuss some simultaneous numerical reconstruction techniques.
\end{remark}

\section{Iterative reconstruction methods}\label{sec:iterative}
We return to the general case in which both $c$ and $\gamma$ are unknown 
\inarxiv{
and consider the Helmholtz model with absorbing boundary conditions. The PML setting could be treated analogously}.

\paragraph{Iteratively regularized Gauss-Newton method IRGNM}
A regularized Gauss-Newton step for solving \eqref{Fxy} defines $\kappagamma^{(n+1)}$ as a minimizer of 
\[
\begin{aligned}
&\|F(\kappagamma^{(n)})+T^{(n)}(\kappagamma-\kappagamma^{(n)})-y\|_{L^2}^2
\\&+\alpha^{(n)}\|(\kappagamma-\kappagamma^{(0)})\|_{L^2}^2
\end{aligned}
\]
where $T^{(n)}=F'(\kappagamma^{(n)})$
thus, with $^*$ denoting the Hilbert space adjoint in $L^2$, the Newton step reads as 
\[
\begin{aligned}
&\kappagamma^{(n+1)} = \kappagamma^{(n)}+({T^{(n)}}^* T^{(n)}+\alpha^{(n)} I)^{-1} 
\\&\hspace*{0.5cm}\Bigl({T^{(n)}}^* (y-F(\kappagamma^{(n)})+\alpha^{(n)}(\kappagamma^{(0)}-\kappagamma^{(n)})\Bigr)
\end{aligned}
\]
or as $\kappagamma^{(n+1)}=\kappagamma^{(n)}+\delta \kappagamma$ where $\delta \kappagamma$ solves the variational equation
\begin{equation}\label{eq:redIRGNM0}
\begin{aligned}
&\langle T^{(n)} \delta\kappagamma, T^{(n)} \xizeta\rangle + \alpha^{(n)} \langle \delta\kappagamma,\xizeta \rangle \\
&= \langle y-F(\kappagamma^{(n)}, T^{(n)} \xizeta\rangle+ \alpha^{(n)} \langle \kappagamma^{(0)}\!\!\!-\kappagamma^{(n)}\!\!,\xizeta \rangle \\ 
&\mbox{for all } \xizeta\in L^2(\widetilde{\Omega})\times L^2(\widetilde{\Omega})\,,
\end{aligned}
\end{equation}
where $\langle\cdot,\cdot\rangle$ denotes the $L^2$ inner products
\[
\begin{aligned}
&\langle \kappagamma,\xizeta\rangle=\langle \kappa,\xi\rangle_{L^2(\widetilde{\Omega})}+\langle \xi,\zeta\rangle_{L^2(\widetilde{\Omega})},\\
&\langle y,z\rangle=\langle y,z\rangle_{L^2(\Gamma)}\,.
\end{aligned}
\]
 
\paragraph{Levenberg-Marquardt method}
A slightly different version of Newton's method is the Levenberg-Marquardt method defining the new iterate as a minimizer of
\[
\begin{aligned}
&\|F(\kappagamma^{(n)})+T^{(n)}(\kappagamma-\kappagamma^{(n)})-y\|_{L^2}^2
\\&+\alpha^{(n)}\|(\kappagamma-\kappagamma^{(n)})\|_{L^2}^2
\end{aligned}
\]
and thus reads as 
\[
\begin{aligned}
&\kappagamma^{(n+1)} = \kappagamma^{(n)}+({T^{(n)}}^* T^{(n)}+\alpha^{(n)} I)^{-1} \\&\hspace*{4cm} {T^{(n)}}^* (y-F(\kappagamma^{(n)})
\end{aligned}
\]
i.e., as $\kappagamma^{(n+1)}=\kappagamma^{(n)}+\delta \kappagamma$ with
\begin{equation}\label{eq:redLM0}
\begin{aligned}
&\langle T^{(n)} \delta\kappagamma, T^{(n)} \xizeta\rangle + \alpha^{(n)} \langle \delta\kappagamma,\xizeta \rangle \\
&= \langle y-F(\kappagamma^{(n)}, T^{(n)} \xizeta\rangle\\ 
&\mbox{for all } \xizeta\in L^2(\widetilde{\Omega})\times L^2(\Omega)\,.
\end{aligned}
\end{equation}

\inarxiv{
\paragraph{All-at-once Newton type methods}
Alternatively, we may apply (regularized versions of) Newton's method directly to \eqref{AuxCu}, i.e., solve the linear system
\begin{equation}\label{eq:AClin}
\begin{aligned}
&A(u^{(n)},\kappagamma^{(n)})+K\,\delta u + L\,\delta\kappagamma=0\\
&Cu^{(n)}+C\delta u=y
\end{aligned}
\end{equation}
for $(\delta u,\delta \kappagamma)=(\delta \hat{\phi}_1,\delta \hat{\phi}_2, \delta \hat{\psi},\delta \kappa, \delta\gamma)$ and set 
\begin{equation}\label{eq:step}
u^{(n+1)}=u^{(n)}+\delta u,\quad  \kappagamma^{(n+1)}=\kappagamma^{(n)}+\delta \kappagamma\,.
\end{equation}
Solving
\begin{equation}\label{eq:unp12}
Ku^{(n+\frac12)}=Ku^{(n)}-A(u^{(n)},\kappagamma^{(n)})
\end{equation}
for $u^{(n+\frac12)}=u^{(n)}-K^{-1}A(u^{(n)},\kappagamma^{(n)})$
and eliminating $\delta u=-K^{-1}(A(u^{(n)},\kappagamma^{(n)}+L\, \delta\kappagamma)=u^{(n+\frac12)}-u^{(n)}-K^{-1}L\, \delta\kappagamma$ in \eqref{eq:AClin} yields an equation for $\delta\kappagamma$ only
\begin{equation}\label{eq:elim}
CK^{-1}L \, \delta \kappagamma = Cu^{(n+\frac12)}-y
\end{equation}
from which, after solving for $\delta\kappagamma$, the new state can be computed as 
\begin{equation}\label{eq:unp1}
u^{(n+1)}=u^{(n+\frac12)}-K^{-1}L\, \delta\kappagamma
\end{equation}
A regularized version of \eqref{eq:elim} is (after computing $u^{(n+\frac12)}$ from \eqref{eq:unp12}) to determine $\delta\kappagamma$ from the variational equation
\begin{equation}\label{eq:aaoIRGNM}
\begin{aligned}
&\langle CK^{-1}L\, \delta\kappagamma, CK^{-1}L\, \xizeta\rangle + \alpha^{(n)} \langle \delta\kappagamma,\xizeta \rangle \\
&= \langle Cu^{(n+\frac12)}-y, CK^{-1}L\, \xizeta\rangle\\
&\qquad
\left[+ \alpha^{(n)} \langle \kappagamma^{(0)}-\kappagamma^{(n)},\xizeta \rangle \right]\\ 
&\mbox{for all } \xizeta\in L^2(\widetilde{\Omega})\times L^2(\Omega)\,,
\end{aligned}
\end{equation}
which is followed by computation of $u^{(n+1)}$ from \eqref{eq:unp1} and of $\kappagamma^{(n+1)}=\kappagamma^{(n)}+\delta\kappagamma$. \\
In \eqref{eq:aaoIRGNM}, the term in large brackets may be skipped to obtain a Levenberg-Marquardt type version of the iteration.

Also the reduced versions \eqref{eq:redIRGNM0}, \eqref{eq:redLM0} can be rewritten in terms of the operators $A$, $C$, $L$, $K$, applying the Implicit Function Theorem to \eqref{eq:S} which yields $S'(\kappagamma^{(n)})=-K^{-1} L$, hence $T^{(n)}=-CK^{-1} L$ with $K,L$ according to \eqref{eq:KL}.
This leads to almost the same formulation as in \eqref{eq:aaoIRGNM}
\begin{equation}\label{eq:redIRGNM}
\begin{aligned}
&\langle CK^{-1}L\, \delta\kappagamma, CK^{-1}L\, \xizeta\rangle + \alpha^{(n)} \langle \delta\kappagamma,\xizeta \rangle \\
&= \langle Cu^{(n)}-y, CK^{-1}L\, \xizeta\rangle\\
&\qquad\left[+ \alpha^{(n)} \langle \kappagamma^{(0)}-\kappagamma^{(n)},\xizeta \rangle \right]\\ 
&\mbox{for all } \xizeta\in L^2(\widetilde{\Omega})\times L^2(\Omega)\,,
\end{aligned}
\end{equation}
(note that here $u^{(n)}$ is defined as the solution to $A(u^{(n)},\kappagamma^{(n)})=0$ and therefore $u^{(n+\frac12)}=u^{(n)}$).
The key difference between \eqref{eq:aaoIRGNM} and \eqref{eq:redIRGNM} lies in the fact that in \eqref{eq:redIRGNM}, the state $u^{(n)}$ has to be precomputed as a solution to the nonlinear state equation $A(u^{(n)},\kappagamma^{(n)})=0$ for given $\kappagamma^{(n)}$, whereas in \eqref{eq:aaoIRGNM}, $u^{(n)}$ just comes from the previous iterate according to \eqref{eq:step} (with $n$ replaced by $n-1$). Instead, in \eqref{eq:aaoIRGNM}, $u^{(n+\frac12)}$ has to be precomputed from a linear state equation \eqref{eq:unp12}. This would make a considerable difference between the two methods in case of a fully nonlinear model. However here, due to the fact that the model decouples into linear subproblems, the difference in computational effort is insignificant. 
}

\paragraph{Gradient type methods}
A Landweber step for solving \eqref{Fxy} is defined by a gradient descent step for the least squares functional
\[
\|F(\kappagamma^{(n)})-y\|_{L^2}^2
\]
i.e., by
\[
\begin{aligned}
\kappagamma^{(n+1)} &= \kappagamma^{(n)}+\mu {T^{(n)}}^* (y-F(\kappagamma^{(n)})
\end{aligned}
\]
with an appropriately chosen step size $\mu$.
\inarxiv{
We just mention that Landweber iteration could as well be applied to the all-at-once version \eqref{AuxCu} and would completely avoid PDE solutions, thus itself act as an iterative Helmholtz solver (however, probably a very slow one). 
}

Some remarks on the implementation are in order. For details we refer to, e.g., \cite{KNSbook:2008}.

\inarxiv{
\paragraph{Discretization}
In the Newton type methods described in Tables \ref{alg:redNewton}, \ref{alg:aaoNewton}, we use a discretization $\kappa^{(n)}(x)=\sum_{i=1}^I a_i d^i(x)$, $\gamma^{(n)}(x)=\sum_{j=1}^J b_j e^j(x)$ on subspaces of $L^2(\widetilde{\Omega})$, by bases $\{d^1,\ldots,d^I\}$, $\{e^1,\ldots,e^J\}$.
}

\paragraph{Choice of $\alpha^{(n)}$}
The regularization parameter in the Newton type methods may be simply chosen along a geometric sequence $\alpha^{(n)}=c \rho^n$ for some $c>0$, $\rho\in(0,1)$ in case of the IRGNM versions (both reduced and all-at-once).
For the Levenberg-Marquardt method, the choice is somewhat more complicated, namely it has to balance nonlinear and linearized residual in the sense of an inexact Newton method such that 
\[
\underline{\theta}\|F(\kappagamma^{(n)}-y\|\leq \mbox{res}(\alpha^{(n)}) \leq\overline{\theta}\|F(\kappagamma^{(n)}-y\|
\]
for some constants $0<\underline{\theta}<\overline{\theta}<1$, where \\
$\mbox{res}(\alpha)=\|F(\kappagamma^{(n)})+T^{(n)}(\kappagamma^{(n+1)}(\alpha)-\kappagamma^{(n)})-y\|$, cf. \cite{Hanke:1997}.

\paragraph{Stopping rule}
To avoid unbounded propagation of the measurement noise through the iterations, the methods defined above
have to be stopped at an appropriate index $n$. A widely used and well-investigated method for this is the discrepancy principle, which for a given noise level $\delta$ and a safety factor $\tau>1$ defines $n$ as the first index such that
\[
\|F(\kappagamma^{(n)})-y\|\leq\tau\delta\,.
\]

\paragraph{Multiple observations} 
As we have seen in Section~\ref{sec:uniqueness}, unique recovery of even just one of the two coefficients $c$ and $\gamma$ requires boundary measurements for several frequencies - a fact that is evident from a simple dimension count. 
Also the fact that the focal point where the high frequency beams interact is moved through the region of interest should be taken into account by incorporating multiple excitations.
This corresponds to using Neumann conditions $g_k^\ell$ at transducer locations $\Sigma_k^\ell$ for $\ell\in\{1,\ldots,L\}$.
Finally, several receiver array locations $\Gamma^m$, $m\in\{1,\ldots,M\}$, might be used to recover a single pair of $\kappa$ and $\gamma$.
Thus, we actually deal with a set of several model and observation operators $A^\ell$, $\ell\in\{1,\ldots,L\}$, $C^m$, $m\in\{1,\ldots,M\}$ respectively. Labelling the resulting forward operators $F_p=C^m \circ S^\ell$ and data $y_p$ for $p=(m-1)L+\ell$, we can write the inverse problem of reconstructing $\kappagamma$ as a system of operator equations
\begin{equation}\label{Fxy_sys}
F_p(\kappagamma)=y_p \, \quad p\in\{1,\ldots,P=L\cdot M\}
\end{equation}
and apply Kaczmarz type methods as follows:
\begin{itemize}
\item[(a)] parallelly apply one step of an iterative reconstruction method to each of the equations in \eqref{Fxy_sys} and then combine the resulting reconstructions $\kappagamma_p$ in a proper way, e.g.,
\begin{equation}\label{Kacz_par}
\begin{aligned}
\kappagamma_p^{(n+1)}&=\kappagamma^{(n)}\!\!\!+G_p(\kappagamma^{(n)}\!\!\!,F_p,y_p),\ p=1,\ldots P, \\ 
\kappagamma^{(n+1)}&=\frac{1}{P}\sum_{p=1}^P \kappagamma_p^{(n+1)} 
\end{aligned}\end{equation}
\item[(b)] sequentially perform one step of an iterative reconstruction method in a 
cyclically repeated manner 
\begin{equation}\label{Kacz_seq}
\kappagamma^{(n+1)}=\kappagamma^{(n)}+G_p(\kappagamma^{(n)},F_p,y_p) 
\end{equation}
where $p=\text{mod}(n-1,P)+1$, 
(the order in which the indices $p$ are addressed could as well be randomized)
\end{itemize}
In here, $G_p(\kappagamma^{(n)},F_p,y_p)$ is defined by one of the Newton or gradient steps 
\inarxiv{from Tables \ref{alg:redNewton}, \ref{alg:aaoNewton}, \ref{alg:redLW} below}.

\inarxiv{
\paragraph{Algorithms}
For a pseudocode description of the methods discussed above, 
see Tables \ref{alg:redNewton}, \ref{alg:aaoNewton}, \ref{alg:redLW} 
}.

\section{Outlook}
In this paper we have made some first steps towards putting the problem of vibroacoustic imaging into the mathematical framework of inverse problems and regularization. We have presented a model in frequency domain, proven uniqueness of recovery of the spatially varying nonlinearity parameter $\gamma(x)$ from pressure measurements at multiple frequencies, and derived Newton as well as gradient based reconstruction methods.   

Natural next steps are on one hand to refine and implement the devised numerical methods and on the other hand to answer important analytical questions. Among the latter, there is uniqueness of 
simultaneous reconstruction of $c(x)$ and $\gamma(x)$.
To this end, the use of multiple excitation locations (instead of or in addition to multiple frequencies), corresponding to shifting the focus of the interacting high-frequency beams around the region of interest, needs to be further investigated. 
Moreover, a priori information should be taken into account. Indeed, an important special case is the one of piecewise constant coefficients, in which only the shapes of finitely many subdomains and finitely many values of $c$ and $\gamma$ are to be found: Here one would expect uniqueness even from  boundary data at just a few frequencies, resulting from appropriately chosen excitations.

A computational framework for the reconstruction of piecewise constant coefficients could be based on the by now standard approach of alternatingly recovering the support and the value of inclusions in a homogeneous background. For a simultaneous recovery of both support and value, the known advantages of total variation regularization can be made use of. In case of known parameter values, also regularization by bound constraints (using the known values as bounds) is a promising approach \cite{comp_minIP}.

Concerning forward simulation, we point to the fact that the high frequency waves $\phi_1$, $\phi_2$ have a strongly preferred direction of propagation, which can justify the use of a parabolic approximation, cf., e.g., \cite{Tappert1977}. 
Indeed, for efficient numerical simulation a decomposition approach has been devised in \cite{Malcolmetal2007,Malcolmetal2008} that splits the forward problem into a three components: (a) directed high frequency propagation of the two beams described by $\phi_1$, $\phi_2$, 
(b) nonlinear interaction of these at the focal point, and 
(c) undirected low frequency propagation to the measurement array via $\psi$. 
This could also be implemented in our framework; 
the adjoint equations for Landweber iteration would have to be re-derived for this purpose. 

Also the model itself might have to be modified. Besides the use of a parabolic approximation in phase (a), also fractional damping e.g., \cite{CaiChenFangHolm_survey2018,fracPAT} is relevant in ultrasonics.
 
\section*{Acknowledgment}
The author wishes to thank Alison Malcolm, University of New Foundland for fruitful discussions that have in fact been crucial in setting up the model considered in this paper. 
The work was supported by the Austrian Science Fund {\sc fwf}
under the grants P30054 and DOC78.

\inarxiv{\subsection*{Appendix: Algorithms}

The following modules are used as subroutines to the algorithms:

\begin{itemize}
\item
$\hat{\phi}_k=\mbox{solveHelmholtz}\Sigma_k(\omega_k,\kappa,\gamma,\hat{g}_k,f)$: solve 
\[
\begin{aligned}
&-\omega_k^2\kappatil \hat{\phi}_k-\Delta \hat{\phi}_k=f \mbox{ in }\Omega\,, \\ 
&\partial_\nu \hat{\phi}_k=-\imath\sigma_k \hat{\phi}_k \mbox{ on } \partial\Omega\setminus\Sigma_k\,, \quad
\partial_\nu \hat{\phi}_k=\hat{g}_k \mbox{ on } \Sigma_k
\end{aligned}
\]
\item
$\hat{\psi}=\mbox{solveHelmholtz}(\omega_d,\kappa,\gamma,f)$: solve 
\[
\begin{aligned}
&-\omega_d^2\kappatil \hat{\psi}-\Delta \hat{\psi}=f \mbox{ in }\Omega\,, \\ 
&\partial_\nu \hat{\psi}=-\imath\sigma \hat{\psi} \mbox{ on } \partial\Omega
\end{aligned}
\]
\item
$q=\mbox{solveHelmholtz}\Gamma(\omega_d,\kappa,\gamma,h)$
\[
\begin{aligned}
&-\omega_d^2\kappatil q-\Delta q=0 \mbox{ in }\Omega\setminus\Gamma\,, \\ 
&\partial_\nu q=\imath\sigma q \mbox{ on } \partial\Omega\,, \quad
\left[\partial_\nu q\right]=-\imath\omega_d h \mbox{ on }\Gamma 
\end{aligned}
\]
\item
\mbox{$[M,r]=\mbox{setupGram\&res}(\omega_d,\vec{\Phi},\vec{\Psi},y[,\alpha,a,b])$:}\\
for $i=1:I+J$\\
\hspace*{7mm} for $i'=1:I+J$\\
\hspace*{14mm} \mbox{compute 
$M_{i,i'}=\omega_d^2 \Re \left(\int_\Gamma \Psi^i \overline{\Psi^{i'}}\, ds\right)$,
\ }\\
\hspace*{7mm} end\\
\hspace*{7mm} compute 
$r_i=\Re \left(\int_\Gamma (\imath\omega_d\hat{\psi}-y) \overline{\imath\omega_d\Psi^{i}}\, ds\right)$\\
\hspace*{7mm}\qquad\qquad$\left[+\alpha D \left(\begin{array}{c}a^{(0)}-a\\b^{(0)}-b\end{array}\right)\right]
$\\
end\\
\end{itemize}
The terms in brackets can be skipped for a Levenberg-Marquardt version.

Moreover, in case of the Newton type methods, for given basis functions $d^i$, $e^i$ used in the representation of $\kappa$ and $\gamma$, respectively, we assume to have precomputed the matrix $D$ with entries
\[
D_{i,i'}=\begin{cases}
\Re \left(\int_{\widetilde{\Omega}} d^i \overline{d^{i'}}\, dx\right) \mbox{ if }i,i'\leq I\\
\Re \left(\int_{\Omega} e^{i-I} \overline{e^{i'-I}}\, dx\right) \mbox{ if }i,i'> I \end{cases}
\]

\newpage
\pagestyle{empty}
{\small \begin{table}[ht]
\hrule
\vskip 2mm
\centerline{regularized Gauss Newton / Levenberg Marquardt step:}
\vskip 2mm
\centerline{reduced version}
\vskip 2mm
\hrule
\vskip 2mm
\noindent
{\bf Input:} 
$\kappa^{(n)}$, $\gamma^{(n)}$ (i.e., coefficients $a^{(n)}$, $b^{(n)}$)\\
{\bf Output:}
$\kappa^{(n+1)}$, $\gamma^{(n+1)}$ (i.e., coefficients $a^{(n+1)}$, $b^{(n+1)}$)
\noindent
\vskip 2mm
\hrule
\vskip 2mm
\noindent
\% computation of $(\hat{\phi}_1^{(n)},\hat{\phi}_2^{(n)},\hat{\psi}^{(n)})=S(\kappagamma^{(n)})$:\\
$\hat{\phi}_k^{(n)} = \mbox{solveHelmholtz}\Sigma_k(\omega_k,\kappa^{(n)},\gamma^{(n)},\hat{g}_k,0)$, $k\in\{1,2\}$.\\
$f= 2(\omega_1-\omega_2) \imath\left(
\kappatil^{(n)}\nabla\hat{\phi}_1^{(n)}\cdot\overline{\nabla\hat{\phi}_2^{(n)}} 
+\omega_1\omega_2 \gammatil^{(n)}\hat{\phi}_1^{(n)}\overline{\hat{\phi}_2^{(n)}}\right)$\\
$\hat{\psi}^{(n)}=\mbox{solveHelmholtz}(\omega_1-\omega_2,\kappa^{(n)},\gamma^{(n)},f)$\\
for $i=1:I$\\
\hspace*{7mm}\% computation of $(\Phi_1^i,\Phi_2^i,\Psi^i)=K^{-1}Ld^i$:\\
\hspace*{7mm} $f^i_k=-\omega_k^2 \hat{\phi}_k^{(n)}\, \chi_{\widetilde{\Omega}}d^i$, \\
\hspace*{7mm} \mbox{$\Phi^i_k= \mbox{solveHelmholtz}\Sigma_k(\omega_k,\kappa^{(n)},\gamma^{(n)},0,f^i_k)$, $k\in\{1,2\}$}\\
\hspace*{7mm} $f^i=-\Bigl((\omega_1-\omega_2)^2 \hat{\psi}
+2(\omega_1-\omega_2) \imath \nabla\hat{\phi}_1^{(n)}\cdot\overline{\nabla\hat{\phi}_2^{(n)}}\Bigr) \chi_{\widetilde{\Omega}}\, d^i$\\
\hspace*{7mm} \mbox{$+2(\omega_1-\omega_2) \imath\Bigl(
\kappatil^{(n)}
\Bigl(\nabla\Phi^i_1\cdot\overline{\nabla\hat{\phi}_2^{(n)}}
+\nabla\hat{\phi}_1^{(n)}\cdot\overline{\nabla\Phi^i_2}\Bigr)$}\\
\hspace*{7mm} \mbox{$+\omega_1\omega_2 \gammatil^{(n)}
\Bigl(\Phi^i_1 \overline{\hat{\phi}_2^{(n)}}
+\hat{\phi}_1^{(n)} \overline{\Phi^i_2}\Bigr)
\Bigr)$}
\\
\hspace*{7mm} $\Psi^i=\mbox{solveHelmholtz}(\omega_1-\omega_2,\kappa^{(n)},\gamma^{(n)},f^i)$\\
end\\
for $j=1:J$ 
\\
\hspace*{7mm}\% computation of $(0,0,\Psi^{I+j})=K^{-1}Le^j$:\\
\hspace*{7mm} $f^j=
-2(\omega_1-\omega_2) \imath \omega_1\omega_2 \, \chi_{\widetilde{\Omega}}e^j \, \hat{\phi}_1^{(n)}\overline{\hat{\phi}_2^{(n)}}
$\\
\hspace*{7mm} $\Psi^{I+j}=\mbox{solveHelmholtz}(\omega_1-\omega_2,\kappa^{(n)},\gamma^{(n)},f^j)$.\\
end\\
\mbox{$[M,r]=\mbox{setupGram\&res}(\omega_1-\omega_2,\vec{\Phi},\vec{\Psi},y[,\alpha^{(n)},a^{(n)},b^{(n)}])$}\\
solve $(M+\alpha^{(n)}D)x=r$,\\
set $a^{(n+1)}=a^{(n)}+(x_1,\ldots,x_I)$, \\
\hspace*{0.4cm} $b^{(n+1)}=b^{(n)}+(x_{I+1},\ldots,x_{I+J})$\\
set $\kappa^{(n+1)}(x)=\sum_{i=1}^I a_i^{(n+1)} d^i(x)$, \\ 
\hspace*{0.5cm}$\gamma^{(n+1)}(x)=\sum_{j=1}^J b_j^{(n+1)} e^j(x)$
\caption{\label{alg:redNewton}}
\end{table}
}

\newpage
\pagestyle{empty}
{\small \begin{table}[ht]
\hrule
\vskip 2mm
\centerline{regularized Gauss Newton / Levenberg Marquardt step: }
\centerline{all-at-once version}
\vskip 2mm
\hrule
\vskip 2mm
\noindent
{\bf Input:} 
\mbox{$\hat{\phi}_1^{(n)}$, $\hat{\phi}_2^{(n)}$, $\hat{\psi}^{(n)}$, $\kappa^{(n)}$, $\gamma^{(n)}$ (i.e., coefficients $a^{(n)}$, $b^{(n)}$)}\\
{\bf Output:}
$\hat{\phi}_1^{(n+1)}$, $\hat{\phi}_2^{(n+1)}$, $\hat{\psi}^{(n+1)}$, $\kappa^{(n+1)}$, $\gamma^{(n+1)}$ \\
(i.e., coefficients $a^{(n+1)}$, $b^{(n+1)}$)
\noindent
\vskip 2mm
\hrule
\vskip 2mm
\noindent
\% computation of $(\hat{\phi}_1^{(n+\frac12)},\hat{\phi}_2^{(n+\frac12)},\hat{\psi}^{(n+\frac12)})=K^{-1}\Bigl(K(\hat{\phi}_1^{(n)},\hat{\phi}_2^{(n)},\hat{\psi}^{(n)})-A(\hat{\phi}_1^{(n)},\hat{\phi}_2^{(n)},\hat{\psi}^{(n)},\kappagamma^{(n)})\Bigr)$:\\
\mbox{$\hat{\phi}_k^{(n+\frac12)}=\mbox{solveHelmholtz}\Sigma_k(\omega_k,\kappa^{(n)},\gamma^{(n)},-\hat{g}_k,0)$, $k\in\{1,2\}$}\\
$f=2(\omega_1-\omega_2) \imath$\\
\mbox{$\cdot\Bigl(
\kappatil^{(n)}\Bigl(\nabla\hat{\phi}_1^{(n+\frac12)}\cdot\overline{\nabla\hat{\phi}_2^{(n)}}+\nabla\hat{\phi}_1^{(n)}\cdot\overline{\nabla\hat{\phi}_2^{(n+\frac12)}}-\nabla\hat{\phi}_1^{(n)}\cdot\overline{\nabla\hat{\phi}_2^{(n)}}\Bigr)$}\\
\hspace*{5mm}$+\omega_1\omega_2 \gammatil^{(n)}\Bigl(\hat{\phi}_1^{(n+\frac12)}\overline{\hat{\phi}_2^{(n)}}+\hat{\phi}_1^{(n)}\overline{\hat{\phi}_2^{(n+\frac12)}}-\hat{\phi}_1^{(n)}\overline{\hat{\phi}_2^{(n)}}\Bigr)\Bigr)$\\  
\mbox{$\hat{\psi}^{(n+\frac12)}=\mbox{solveHelmholtz}(\omega_1-\omega_2,\kappa^{(n)},\gamma^{(n)},f)$}\\
for $i=1:I$\\
\hspace*{7mm}\% computation of $(\Phi_1^i,\Phi_2^i,\Psi^i)=K^{-1}Ld^i$:\\
\hspace*{7mm} $f^i_k=-\omega_k^2 \hat{\phi}_k^{(n)}\, \chi_{\widetilde{\Omega}}d^i$, \\
\hspace*{7mm} \mbox{$\Phi^i_k=\mbox{solveHelmholtz}\Sigma_k(\omega_k,\kappa^{(n)},\gamma^{(n)},0,f^i_k)$, $k\in\{1,2\}$.}\\
\hspace*{7mm} $f^i=-\Bigl((\omega_1-\omega_2)^2 \hat{\psi}
+2(\omega_1-\omega_2) \imath \nabla\hat{\phi}_1^{(n)}\cdot\overline{\nabla\hat{\phi}_2^{(n)}}\Bigr) \chi_{\widetilde{\Omega}}\, d^i$\\
\hspace*{7mm} \mbox{$+2(\omega_1-\omega_2) \imath\Bigl(
\kappatil^{(n)}
\Bigl(\nabla\Phi^i_1\cdot\overline{\nabla\hat{\phi}_2^{(n)}}
+\nabla\hat{\phi}_1^{(n)}\cdot\overline{\nabla\Phi^i_2}\Bigr)$}
\\
\hspace*{7mm} \qquad\qquad\mbox{$+\omega_1\omega_2 \gammatil^{(n)}
\Bigl(\Phi^i_1 \overline{\hat{\phi}_2^{(n)}}
+\hat{\phi}_1^{(n)} \overline{\Phi^i_2}\Bigr)
\Bigr)$}
\\
\hspace*{7mm} \mbox{$\Psi^i=\mbox{solveHelmholtz}(\omega_1-\omega_2,\kappa^{(n)},\gamma^{(n)},f^i)$}\\
end\\
for $j=1:J$ \%(note that $\Phi^{I+j}=0$)
\\
\hspace*{7mm}\% computation of $(0,0,\Psi^{I+j})=K^{-1}Le^j$:\\
\hspace*{7mm} $f^j=
-2(\omega_1-\omega_2) \imath \omega_1\omega_2 \, \chi_{\widetilde{\Omega}}e^j \, \hat{\phi}_1^{(n)}\overline{\hat{\phi}_2^{(n)}}
$\\
\hspace*{7mm} $\Psi^{I+j}=\mbox{solveHelmholtz}(\omega_1-\omega_2,\kappa^{(n)},\gamma^{(n)},f^j)$.\\
end\\
\mbox{$[M,r]=\mbox{setupGram\&res}(\omega_1-\omega_2,\vec{\Phi},\vec{\Psi},y[,\alpha^{(n)},a^{(n)},b^{(n)}])$}\\
solve $(M+\alpha^{(n)}D)x=r$,\\
set $a^{(n+1)}=a^{(n)}+(x_1,\ldots,x_I)$, \\
\hspace*{0.4cm} $b^{(n+1)}=b^{(n)}+(x_{I+1},\ldots,x_{I+J})$\\
set $\kappa^{(n+1)}(x)=\sum_{i=1}^I a_i^{(n+1)} d^i(x)$, \\ 
\hspace*{0.5cm}$\gamma^{(n+1)}(x)=\sum_{j=1}^J b_j^{(n+1)} e^j(x)$
\\
\hspace*{0.5cm}$\hat{\phi}_k^{(n+1)}=\hat{\phi}_k^{(n+\frac12)}-\sum_{i=1}^I a_i^{(n+1)} \Phi^i$, $k\in\{1,2\}$,\\
\hspace*{0.5cm}\mbox{$\hat{\psi}^{(n+1)}=\hat{\psi}^{(n+\frac12)}-\sum_{i=1}^I a_i^{(n+1)} \Psi^i-\sum_{j=1}^J b_j^{(n+1)} \Psi^{I+j}$}
\caption{\label{alg:aaoNewton}}
\end{table}
}

\newpage
\pagestyle{empty}
{\small \begin{table}[ht]
\hrule
\vskip 2mm
\centerline{reduced Landweber step}
\vskip 2mm
\hrule
\vskip 2mm
\noindent
{\bf Input:} 
$\kappa^{(n)}$, $\gamma^{(n)}$ \\
{\bf Output:}
$\kappa^{(n+1)}$, $\gamma^{(n+1)}$ 
\noindent
\vskip 2mm
\hrule
\vskip 2mm
\noindent
\% computation of $(\hat{\phi}_1^{(n)},\hat{\phi}_2^{(n)},\hat{\psi}^{(n)})=S(\kappagamma^{(n)})$:\\
$\hat{\phi}_k^{(n)} = \mbox{solveHelmholtz}\Sigma_k(\omega_k,\kappa^{(n)},\gamma^{(n)},\hat{g}_k,0)$, $k\in\{1,2\}$.\\
$f= 2(\omega_1-\omega_2) \imath\left(
\kappatil^{(n)}\nabla\hat{\phi}_1^{(n)}\cdot\overline{\nabla\hat{\phi}_2^{(n)}} 
+\omega_1\omega_2 \gammatil^{(n)}\hat{\phi}_1^{(n)}\overline{\hat{\phi}_2^{(n)}}\right)$\\
$\hat{\psi}^{(n)}=\mbox{solveHelmholtz}(\omega_1-\omega_2,\kappa^{(n)},\gamma^{(n)},f)$\\
\% computation of adjoint states:\\
$r=\mbox{tr}_\Gamma(\imath(\omega_1-\omega_2)\hat{\psi}^{(n)}-y)$\\
$q^{(n)}=\mbox{solveHelmholtz}\Gamma(\omega_1-\omega_2,\kappa^{(n)},\gamma^{(n)},r)$\\
$f_k=2(\omega_1-\omega_2)\imath\Bigl(
\nabla \bigl( \kappatil^{(n)} \overline{q}\nabla \phi_{\not{k}}^{(n)}\bigr)
+\omega_1\omega_2 \gammatil^{(n)} \overline{q}\phi_{\not{k}}^{(n)}\Bigr)$\\
$p_k^{(n)}=\mbox{solveHelmholtz}\Sigma_k(\omega_k,\kappa^{(n)},\gamma^{(n)},0,f_k)$, $k=1,2$\\
compute
\[
\begin{aligned}
&\xi^{(n)}=-\Re\Bigl((1-\imath)\Bigl(\sum_{k=1}^2 \omega_k^2 \hat{\phi}_k^{(n)}\overline{p_k^{(n)}}
+(\omega_1-\omega_2)^2 \hat{\psi}^{(n)}\overline{q^{(n)}}\\
&\qquad\qquad\qquad+2(\omega_1-\omega_2) \imath \nabla\hat{\phi}_1^{(n)}\cdot\overline{\nabla\hat{\phi}_2^{(n)}} \overline{q^{(n)}} \Bigr)\Bigr)\vert_{\widetilde{\Omega}}
\\
&\zeta^{(n)}=-\Re\Bigl((1-\imath)\Bigl(2(\omega_1-\omega_2) \imath \omega_1\omega_2 \hat{\phi}_1^{(n)}\overline{\hat{\phi}_2^{(n)}}\overline{q^{(n)}}\Bigr)\Bigr)\vert_{\widetilde{\Omega}}
\end{aligned}
\]
set $\kappa^{(n+1)}=\kappa^{(n)}-\mu\xi^{(n)}$, $\gamma^{(n+1)}=\gamma^{(n)}-\mu\zeta^{(n)}$
\caption{\label{alg:redLW}}
\end{table}
}

}



\end{document}